\documentclass{amsart}
\usepackage{latexsym}

\usepackage{amsfonts,amssymb,amsmath}
\usepackage{amsthm}
\usepackage{graphicx,color}
\usepackage{color}
\usepackage{amsmath} 
\newtheorem{theorem}{Theorem}	
\newtheorem{lemma}{Lemma}[section]
		
\newtheorem{proposition}{Proposition}		
\newtheorem{claim}{Claim}[section]	
			
\newtheorem{definition}{Definition}

\newfont{\bg}{cmr9 scaled\magstep2}
\newcommand{\bigzerol}{\smash{%\lower1.0ex
\hbox{\bg 0}}}
\newcommand{\bigzerou}{
\smash{\hbox{\bg 0}}}

\title[Generalized distance-squared mappings of $\mathbb{R}^{n+1}$ into $\mathbb{R}^{2n+1}$]%
{
Generalized distance-squared mappings \\ 
of $\mathbb{R}^{n+1}$ into $\mathbb{R}^{2n+1}$ 
}
\author{S.~Ichiki}
\address{Dai Nippon Printing Co., Ltd., 
Tokyo 162-8001, JAPAN} 
\email{ichiki-shunsuke-jb@ynu.jp}
\author{T.~Nishimura%\footnote{Corresponding author.   }}
} 
\address{
Research Institute of Environment and Information Sciences,  
Yokohama National University, 
Yokohama 240-8501, JAPAN}
\email{nishimura-takashi-yx@ynu.jp}
\begin{document}
\date{}
%%%%%%%%%%%%%%%%%%%%%%%%%%%%%%%%%%%%%%%%%%%%% 
\begin{abstract}
We classify generalized distance-squared mappings 
of $\mathbb{R}^{n+1}$ into $\mathbb{R}^{2n+1}$ ($n\ge 1$) 
having generic central points.      
Moreover, 
we show that there does not exist a universal bad set $\Sigma\subset (\mathbb{R}^{n+1})^{2n+1}$ 
in the case of this dimension-pair.        
\end{abstract}
%%%%%%%%%%%%%%%%%%%%%%%%%%%%%%%%%%%%%%%%%%%%%%  
\maketitle
\noindent
\section{Introduction}\label{Introduction}
%%%%%%%%%%%%%%%%%%%%%%%%%%%%%%%%%%%%%%%%%%%%%%%%%% 
%%%%%%%%%%%%%%%%%%%%%%%%%%%%%%%%%%%%%%%%%%%%%%%%%% 
%Let $n$ be a positive integer.   For any positive integer $k$, let 
For any positive integers $k, n$, 
%Let $p_0, p_1, \ldots, p_k$ be $(k+1)$ points of $\mathbb{R}^{n+1}$. 
let $p_i=(p_{i0}, p_{i1}, \ldots, p_{in})$  $(0\le i\le k)$ 
(resp., $A=(a_{ij})_{0\le i\le k, 0\le j\le n}$) 
be a point of $\mathbb{R}^{n+1}$ 
(resp., a $(k+1)\times (n+1)$ matrix with non-zero entries).     
%Let $A=(a_{ij})_{0\le i\le k, 0\le j\le n}$ be a $(k+1)\times (n+1)$ matrix with non-zero entries.   
Let $G_{(p_0, p_1, \ldots, p_k, A)}: \mathbb{R}^{n+1}\to \mathbb{R}^{k+1}$ be the mapping 
defined by 
%a {\it generalized distance-squared mapping}:   
{\small 
\[
G_{(p_0, p_1, \ldots, p_k, A)}(x)=\left(
\sum_{j=0}^n a_{0j}(x_j-p_{0j})^2, 
\sum_{j=0}^n a_{1j}(x_j-p_{1j})^2, 
\ldots, 
\sum_{j=0}^n a_{kj}(x_j-p_{kj})^2
\right), 
\] 
}
where $x=(x_0, x_1, \ldots, x_n)\in \mathbb{R}^{n+1}$.    
The mapping $G_{(p_0, p_1, \ldots, p_k, A)}$ is called a {\it generalized distance-squared mapping}.    
Each component of a generalized distance-squared mapping defines the family of quadrics.   
The singularities of  $G_{(p_0, p_1, \ldots, p_k, A)}$ is a helpful information on the contacts of these families.   
Thus, we may regard that generalized distance-squared mappings are a significant tool 
in the applications of singularity theory to differential geometry.    
Therefore, it is natural to classify them.    
%Generalized distance-squared mappings are a useful tool in the
%applications of singularity theory to differential geometry. Their
%singularities give information on the contacts amongst the families
%of quadrics defined by the components of
%$G_{(p_0,p_1,\ldots,p_k,A)}$.     
%It is therefore natural to classify
%mappings $G_{(p_0,p_1,\ldots,p_k,A)}$.    
In the case of $n=k=1$, a recognizable classification is known for $G_{(p_0,p_1,\ldots,p_k,A)}$ 
(see Proposition \ref{proposition 1} below).   
The purpose of this paper is to obtain a classification in the case of $k=2n$.     
\par  
Two mappings $f_i: \mathbb{R}^{n+1}\to \mathbb{R}^{k+1}$ $(i=1,2)$
%(resp., map-germs $f_i: (\mathbb{R}^{n+1}, q_i)\to (\mathbb{R}^{k+1},f_i(q_i)$)  $(i=1, 2)$ 
are said to be 
{\it $\mathcal{A}$-equivalent} if there exist $C^\infty$ diffeomorphisms 
$h: \mathbb{R}^{n+1}\to \mathbb{R}^{n+1}$ and $H:\mathbb{R}^{k+1}\to \mathbb{R}^{k+1} $
%(resp.,  germs of $C^\infty$ diffeomorphism 
%$h: (\mathbb{R}^{n+1}, q_1)\to (\mathbb{R}^{n+1},q_2)$ and 
%$H: (\mathbb{R}^{k+1}, f_2(q_2))\to (\mathbb{R}^{k+1},f_1(q_1))$) 
such that 
the identity $f_1=H\circ f_2\circ h$ holds.     
%%%%%%%%%%%%%%%%%%%%%%%%%%%%%%%%%%%%%%%%%%%%
%%%%%%%%%%%%%%%%%%%%%%%%%%%%%%%%%%%%%%%%%%%% 
\begin{proposition}[\cite{generalized}] \label{proposition 1}
Let $((x_0, y_0), (x_1, y_1))$ be the standard coordinates of $(\mathbb{R}^2)^2$ 
and let $\Sigma$ be the hypersurface in $(\mathbb{R}^2)^2$ defined by 
$(x_0-x_1)(y_0-y_1)=0$.    
Let $(p_0, p_1)$ be a point in $(\mathbb{R}^2)^2 - \Sigma$ and 
let $A_k$ be a $2\times 2$ matrix of rank $k$ with non-zero entries (k=1, 2).      
Then, the following hold:   
\begin{enumerate}
\item 
%Suppose that the rank of $A$ is $1$.   
%Then, 
The mapping $G_{(p_0, p_1, A_1)}$ is proper and stable, and it is 
$\mathcal{A}$-equivalent to $(x,y)\mapsto (x,y^2)$.   
%\item 
%Suppose that the rank of $A$ is $2$.   
\item 
The mapping $G_{(p_0, p_1, A_2)}$ is proper and stable, and 
it is not $\mathcal{A}$-equivalent to  $G_{(p_0, p_1, A_1)}$.   
%with only one cusp point having two connected components of singular set.   
\item Let $B_2$ be a $2\times 2$ matrix of rank $2$ with non-zero entries and let 
$(q_0, q_1)$ be a point in $(\mathbb{R}^2)^2 - \Sigma$.   
Then, $G_{(p_0, p_1, A_2)}$ is $\mathcal{A}$-equivalent 
to $G_{(q_0, q_1, B_2)}$. 
%but not $\mathcal{A}$-equivalent to $\Phi_2$, 
%where $\widetilde{p}_0=(0,0), \widetilde{p}_0=(1,1)$ and 
%$A_1=\left(
%\begin{array}{rr}
% 1 & 1 \\ 
%-1 & 1 
%\end{array}
%\right).
%$   
\end{enumerate}
\end{proposition}
%%%%%%%%%%%%%%%%%%%%%%%%%%%%%%%%%%%%%%%%%%%%%%%%%%%%%%%%%%%%%%  
\par 
In the case of $k=2n$,  
only a few partial classification result for $G_{(p_0, \ldots, p_k, A)}$ is known.    
A {\it distance-squared mapping}   
$D_{(p_0, p_1, \ldots, p_k)}$ 
(resp., {\it Lorentzian distance-squared mappings}  
$L_{(p_0, p_1, \ldots, p_k)}$) 
is the mapping 
$G_{(p_0, p_1, \ldots, p_k,A)}$ satisfying that each entry of $A$ is $1$ 
(resp., $a_{i0}=-1 $ and $a_{ij}=1$ if $j\ne 0$).   
%%%%%%%%%%%%%%%%%%%%%%%%%%%%%%%%%%%%%%%%%%%%%%%%%%%%   
\begin{proposition}[\cite{ichikinishimura, ichikinishimura2}]
\label{proposition 2}
There exists a closed subset $\Sigma\subset (\mathbb{R}^{n+1})^{2n+1}$ 
with Lebesgue measure zero 
such that for any $p=(p_0, p_1, \ldots, p_{2n})\in (\mathbb{R}^{n+1})^{2n+1}-\Sigma$,  
both $D_{(p_0, p_1, \ldots, p_{2n})}$ and $L_{(p_0, p_1, \ldots, p_{2n})}$ 
are $\mathcal{A}$-equivalent to the inclusion  
$(x_0, x_1, \ldots, x_n)\mapsto (x_0, x_1, \ldots, x_n, 0, \ldots, 0).$   
\end{proposition}
%%%%%%%%%%%%%%%%%%%%%%%%%%%%%%%%%%%%%%%%%%%%%%%%%%% 

%The following Theorem \ref{theorem 1} generalizes Proposition \ref{proposition 2}.    
\begin{theorem} 
\label{theorem 1}
Let $A=(a_{ij})_{0\le i\le 2n, 0\le j\le n}$ be a $(2n+1)\times (n+1)$ matrix with non-zero entries.   
Then, the following two hold: 
\begin{enumerate}
\item Suppose that the rank of $A$ is $n+1$.    Then, there exists a closed subset $\Sigma_A\subset (\mathbb{R}^{n+1})^{2n+1}$ 
with Lebesgue measure zero 
such that for any $p=(p_0, p_1, \ldots, p_{2n})\in  (\mathbb{R}^{n+1})^{2n+1}-\Sigma_A$, 
$G_{(p, A)}$ is $\mathcal{A}$-equivalent to the following mapping:   
\[
(x_0, x_1, \ldots, x_n)\mapsto 
(x_0^2, x_0x_1, \ldots, x_0x_n, x_1, \ldots, x_n).   
\]
\item Suppose that the rank of $A$ is less than $n+1$.    
Then, there exists a closed subset $\Sigma_A\subset (\mathbb{R}^{n+1})^{2n+1}$ with Lebesgue measure zero 
such that for any $p=(p_0, p_1, \ldots, p_{2n})\in  (\mathbb{R}^{n+1})^{2n+1}-\Sigma_A$, 
$G_{(p, A)}$ is $\mathcal{A}$-equivalent to the inclusion 
$(x_0, x_1, \ldots, x_n)\mapsto (x_0, x_1,\ldots, x_n, 0, \ldots, 0).$
\end{enumerate}
\end{theorem} 
The mapping given in the assertion (1) of Theorem \ref{theorem 1} was firstly given in 
\cite{whitney} and is called {\it the normal form of Whitney umbrella}.      
%defined in \cite{whitney}.     
It is easily seen that the normal form of Whitney umbrella is not $\mathcal{A}$-equivalent to the inclusion 
$(x_0, x_1, \ldots, x_n)\mapsto (x_0, x_1, \ldots, x_n, 0, \ldots, 0)$.    
Moreover, by Mather's characterization 
theorem of stable mappings given in \cite{mather5}, 
it is easily shown that these two mappings are proper and stable.         
Thus, Theorem \ref{theorem 1} may be regarded as a result of Proposition \ref{proposition 1} type.    
On the other hand, it is desirable to improve Theorem \ref{theorem 1} so that the bad set $\Sigma_A$ given in Theorem \ref{theorem 1} does not depend on the given matrix $A$.  
However, contrary to the case of $n=k=1$, 
in this case it is impossible to expect the existence of such a universal bad set $\Sigma$ as follows.   
\begin{theorem} 
\label{theorem 2}
%Let $A=(a_{ij})_{0\le i\le 2n, 0\le j\le n}$ be a $(2n+1)\times (n+1)$ matrix with positive entries.   
There does not  exist a closed subset $\Sigma\subset (\mathbb{R}^{n+1})^{2n+1}$  with Lebesgue measure zero 
such that for any point $p=(p_0, p_1, \ldots, p_{2n})\in  (\mathbb{R}^{n+1})^{2n+1}-\Sigma$, 
the following two hold.    
\begin{enumerate}
\item Suppose that $A$ is a 
$(2n+1)\times (n+1)$ matrix with non-zero entries such that 
the rank of $A$ is $n+1$.    Then, 
%there exists a closed subset $\Sigma_A\subset (\mathbb{R}^{n+1})^{2n+1}$ 
%with Lebesgue measure zero 
%such that for any $p=(p_0, \ldots, p_{2n+1})\in  (\mathbb{R}^{n+1})^{2n+1}-\Sigma_A$, 
$G_{(p, A)}$ is $\mathcal{A}$-equivalent to the following mapping:   
\[
(x_0, x_1, \ldots, x_n)\mapsto 
(x_0^2, x_0x_1, \ldots, x_0x_n, x_1, \ldots, x_n).   
\]
\item Suppose that  $A$ is a 
$(2n+1)\times (n+1)$ matrix with non-zero entries such that 
the rank of $A$ is less than $n+1$.    
Then, 
%there exists a closed subset $\Sigma_A\subset (\mathbb{R}^{n+1})^{2n+1}$ with Lebesgue measure zero 
%such that for any $p=(p_0, \ldots, p_{2n+1})\in  (\mathbb{R}^{n+1})^{2n+1}-\Sigma_A$, 
$G_{(p, A)}$ is $\mathcal{A}$-equivalent to the inclusion 
$(x_0, x_1, \ldots, x_n)\mapsto (x_0, x_1, \ldots, x_n, 0, \ldots, 0).$
\end{enumerate}
\end{theorem} 
\par 
\bigskip 
%%%%%%%%%%%  
The assertion (1) of Theorem \ref{theorem 1}, the assertion (2) of Theorem \ref{theorem 1} 
and Theorem \ref{theorem 2} 
%, and \ref{theorem 2} 
are 
%given.   
%Theorems \ref{theorem 1} and \ref{theorem 2} 
proved in Sections \ref{section 2}, \ref{section 3}
and \ref{section 4} 
%, and \ref{section 5} 
respectively.   
%%%%%%%%%%%%%%%%%%%%%%%%%%%%%%%%%%%%%%%%%%%%%%%%%% 
%%%%%%%%%%%%%%%%%%%%%%%%%%%%%%%%%%%%%%%%%%%%%%%%%%    
\section{Proof of the assertion (1) of Theorem \ref{theorem 1}}\label{section 2}
%%%%%%%%%%%%%%%%%%%%%%%%%%%%%%%%%%%%%%%%%%%%%%%%%% 
Set   ${A}_1=(a_{ij})_{0\le i, j\le n }$ and ${A}_2=(a_{ij})_{n+1\le i\le 2n, 0\le j\le n}$.     
Taking permutations of coordinates of the target space if necessary,  without loss of generality, 
from the first we may assume that rank$(A_1)=n+1$. 
%%%%%%%%%%%%%%%%%%%%%%%%%%%%%%%%%%%%%%%%%%%%%%%%%%    
\subsection{STEP 1}\label{subsection 2.1}
%%%%%%%%%%%%%%%%%%%%%%%%%%%%%%%%%%%%%%%%%%%%%%%%%%
The purpose of this step is to delete quadratic terms as many as possible 
by a linear transformation $H_1: \mathbb{R}^{2n+1}\to \mathbb{R}^{2n+1}$ of the following type.    
\begin{eqnarray*}
{ } & { } & H_1(X_0, \ldots, X_{2n}) \\ 
{ } & = & (X_0, \ldots, X_{2n}) 
\left(
\begin{array}{cccccc}
\lambda_{0,0} & \cdots & \lambda_{0,n} & \lambda_{0,n+1} & \cdots & \lambda_{0,2n} \\ 
\vdots & \ddots & \vdots & \ddots & \vdots & \vdots \\ 
\lambda_{n,0} & \cdots  & \lambda_{n,n} & \lambda_{n,n+1} & \cdots & \lambda_{n,2n} \\ 
0 & \cdots & 0 & 1 & 0 & 0 \\ 
\vdots & \ddots & \vdots & 0 & \ddots & 0 \\ 
0 & \cdots & 0 & 0 & 0 & 1  
\end{array}
\right).   
\end{eqnarray*}
Set ${\Lambda}_1=(\lambda_{i,j})_{0\le i, j\le n}$ and ${\Lambda}_2=(\lambda_{i,j})_{0\le i\le n, n+1\le j\le 2n}$.     
Two matrices ${\Lambda}_1$ and ${\Lambda}_2$ are obtained as the solutions of the following linear equations, 
where $E_{n+1}$ is the $(n+1)\times (n+1)$ unit matrix and $M^T$ stands for the transposed matrix of a matrix $M$.    
\[
A_1^T{\Lambda}_1=E_{n+1}, \quad 
A_1^T{\Lambda}_2=-A_2^T.   
\]  
Set $H_1\circ G_{(p, A)}=(\varphi_{1,0}, \varphi_{1,1}, \ldots, \varphi_{1,2n})$.    
Then, $\varphi_{1, i}$ $(0\le i\le 2n)$ may be expressed as follows.    
\[
\varphi_{1, i}(x_0, x_1, \ldots, x_n)  
= 
\left\{ 
\begin{array}{ll}
x_i^2+\sum_{j=0}^n b_{ij}x_j+c_i &  (0\le i\le n) \\
\sum_{j=0}^n b_{ij}x_j+c_i &  (n+1\le i\le 2n), 
\end{array} 
\right.
\]
where $c_i$ stands for the constant term and $b_{ij}$ is as follows.   
\[
b_{ij}  
= 
\left\{ 
\begin{array}{ll}
-2\sum_{k=0}^n \lambda_{k,i}a_{kj}{p_{kj}} &  (0\le i\le n) \\
-2\left(\sum_{k=0}^n \lambda_{k,i}a_{kj}p_{kj} + a_{ij}{p_{ij}}\right) &  (n+1\le i\le 2n).   
\end{array} 
\right.
\]
%%%%%%%%%%%%%%%%%%%%%%%%%%%%%%%%%%%%%%%%%%%%%%%%%%%%%%%%%%%%%%%%%%% 
\subsection{STEP 2}\label{subsection 2.2} 
%%%%%%%%%%%%%%%%%%%%%%%%%%%%%%%%%%%%%%%%%%%%%%%%%%%%%%%%%%%%%%%%%%% 
The purpose of this step is to delete constant terms by the parallel transformation $H_2: \mathbb{R}^{2n+1}\to \mathbb{R}^{2n+1}$ 
defined by 
\[
H_2(X_0, X_1, \ldots, X_{2n})=
(X_0-c_0, X_1-c_1, \ldots, X_{2n}-c_{2n}). 
\]    
Set $H_2\circ H_1\circ G_{(p, A)}=(\varphi_{2,0}, \varphi_{2,1}, \ldots, \varphi_{2,2n})$.    
Then, $\varphi_{2, i}$ $(0\le i\le 2n)$ may be expressed as follows.    
\[
\varphi_{2, i}(x_0, x_1, \ldots, x_n)  
= 
\left\{ 
\begin{array}{ll}
x_i^2+\sum_{j=0}^n b_{ij}x_j &  (0\le i\le n) \\
\sum_{j=0}^n b_{ij}x_j &  (n+1\le i\le 2n), 
\end{array} 
\right.
\]
%%%%%%%%%%%%%%%%%%%%%%%%%%%%%%%%%%%%%%%%%%%%%%%%%%    
\subsection{STEP 3}\label{subsection 2.3}
%%%%%%%%%%%%%%%%%%%%%%%%%%%%%%%%%%%%%%%%%%%%%%%%%%
The purpose of this step is to construct the bad set $\Sigma_A$. 
The desirable bad set $\Sigma_A$ has the property that for any 
$p\in (\mathbb{R}^{n+1})^{2n+1}-\Sigma_A$, 
$H_2\circ H_1\circ G_{(p,A)}$ is transformed to the following type form:  
\[
\left(
x_0^2+d_{0,0}x_0, \cdots, x_n^2+d_{n,n}x_n, 
d_{n+1,0}x_0+d_{n+1, 1}x_1, \cdots, 
d_{2n, 0}x_0+d_{2n, n}x_n
\right). 
\]
by composing a linear transformation $H_3: \mathbb{R}^{2n+1}\to \mathbb{R}^{2n+1}$ of the following type.    
\begin{eqnarray*}
{ } & { } & H_3(X_0, X_1, \ldots, X_{2n}) \\ 
{ } & = & (X_0, X_1, \ldots, X_{2n}) 
\left(
\begin{array}{cccccc}
1 & 0 & 0 & 0 & \cdots & 0 \\ 
0 & \ddots & 0 & \vdots & \vdots & \vdots \\ 
0 & 0  & 1 & 0 & \cdots & 0 \\ 
\gamma_{n+1, 0} & \cdots & \gamma_{n+1, n} & \gamma_{n+1, n+1} & \cdots & \gamma_{n+1,2n} \\ 
\vdots & \vdots & \vdots & \vdots & \vdots & \vdots \\ 
\gamma_{2n, 0} & \cdots & \gamma_{2n, n} & \gamma_{2n, n+1} & \cdots & \gamma_{2n,2n}   
\end{array}
\right), 
\end{eqnarray*}
where $d_{ii}$ $(0\le i \le n)$ are constants and $d_{ij}$ $(n+1\le i \le 2n)$ are non-zero constants.   
\par 
%\bigskip 
For any $j$ $(0\le j\le n)$, set 
${\bf b}_j=(b_{n+1, j}, \ldots, b_{2n, j})^T$.     Moreover, for any  $j$ $(0\le j\le n)$, set 
\[
B_j=({\bf b}_0, \ldots, {\bf b}_{j-1}, \hat{{\bf b}_j}, {\bf b}_{j+1}, \ldots, {\bf b}_n), 
\] 
where  $\hat{{\bf b}_j}$ stands for deleting ${\bf b}_j$.    
Then, $B_j$ is an $n\times n$ matrix for any $j$ $(0\le j\le n)$.     
\begin{definition}
{\rm 
\begin{enumerate}
\item For any  $j$ $(0\le j\le n)$, $\Sigma_{B_j}$ is the set consisting of points $p\in (\mathbb{R}^{n+1})^{2n+1}$ 
such that $\det B_j=0$.     
\item 
\[
\Sigma_A=\bigcup_{j=0}^n \Sigma_{B_j}.   
\]
\end{enumerate}
}
\end{definition}
The set $\Sigma_A$ is closed and 
of Lebesgue measure zero since it is an algebraic set. 
% in $(\mathbb{R}^{n+1})^{2n+1}$.   
Set $H_3\circ H_2\circ H_1\circ G_{(p, A)}=(\varphi_{3,0}, \varphi_{3,1}, \ldots, \varphi_{3,2n})$.    
\begin{lemma}
\label{lemma 2.1}
Let $p$ be a point of $(\mathbb{R}^{n+1})^{2n+1}-\Sigma_A$.    
Then, the following hold:   
\begin{enumerate}
\item 
For any $j$ $(0\le j\le n)$, 
$\varphi_{3, j}$ may be expressed as follows where $d_{j, j}$ is a constant.    
\[
\varphi_{3, j}(x_0, x_1, \ldots, x_n)
= 
x_j^2+d_{j,j}x_j.
\]
\item 
For any $j$ $(n+1\le j\le 2n)$, 
$\varphi_{3, j}$ may be expressed as follows where $d_{j,0}, d_{j, j-n}$ are non-zero constants. 
\[
\varphi_{3, j}(x_0, x_1, \ldots, x_n)
= 
 d_{j, 0}x_0+d_{j, j-n}x_{j-n}. 
\] 
\end{enumerate}
\end{lemma}
\noindent 
\proof\quad 
Since the given point $p$ is outside the bad set $\Sigma_A$, it follows that 
$\det B_j\ne 0$ for any $j$ $(0\le j\le n)$.    
Thus, for any $j$ $(0\le j\le n)$, there exists the unique solution 
$(\gamma_{n+1, j}, \ldots, \gamma_{2n, j})^T$ for the following linear equations.  
\[
(-b_{j, 0}, \ldots, -b_{j, j-1}, \hat{-b_{j, j}}, -b_{j, j+1}, \ldots, -b_{j, n})^T= 
B_j^T(\gamma_{n+1, j}, \ldots, \gamma_{2n, j})^T,    
\] 
where $\hat{-b_{j, j}}$ stands for deleting $-b_{j, j}$.   
Therefore, the assertion (1) follows.     
\par 
Next, we show the assertion (2).       
%For any $j$ $(0\le j\le n)$, let $(\gamma_{n+1}, \ldots, \gamma_{2n})^T$ is a non-zero 
%vector which is perpendicular to ${\bf b}_j$.  
For any  $j$ $(1\le j\le n)$, set 
\[
\widetilde{B}_j=({\bf b}_1\ldots, \hat{{\bf b}_j}, \ldots, {\bf b}_n)
\] 
where  $\hat{{\bf b}_j}$ stands for deleting ${\bf b}_j$.    
Then, $\widetilde{B}_j$ is an $n\times (n-1)$ matrix for any $j$ $(1\le j\le n)$.   
Since  $\det B_j\ne 0$ for any $j$ $(0\le j\le n)$, 
it follows that rank$(\widetilde{B}_j)=n-1$ 
for any  $j$ $(1\le j\le n)$.  
Thus, for any  $j$ $(1\le j\le n)$, 
the solution space for the linear equation  
\[
(0, \ldots, 0)^T=\widetilde{B}_j^T(\gamma_{n+1, n+j}, \ldots, \gamma_{2n, n+j})^T
\] 
is one-dimensional.    
Moreover, since  $\det B_j\ne 0$ for any $j$ $(0\le j\le n)$, both 
${\bf b}_0$ and ${\bf b}_j$ are not perpendicular to non-zero 
solution vector  $(\gamma_{n+1, n+j}, \ldots, \gamma_{2n, n+j})$ 
for the above linear equations 
for any $j$ $(1\le j\le n)$.   
Thus, the assertion (2) follows.    
\hfill $\Box$ 
%%%%%%%%%%%%%%%%%%%%%%%%%%%%%%%%%%%%%%%%%%%%%%%%%%%%%%%%%%%%%%%%%%% 
\subsection{STEP 4}\label{subsection 2.4} 
%%%%%%%%%%%%%%%%%%%%%%%%%%%%%%%%%%%%%%%%%%%%%%%%%%%%%%%%%%%%%%%%%%% 
The purpose of this step is to complete the square for the first $(n+1)$ components of 
$H_3\circ H_2\circ H_1\circ G_{(p,A)}$.     
Set 
\[
h_1(x_0, x_1, \ldots, x_n)=(x_0-\frac{1}{2}d_{0,0}, x_1-\frac{1}{2}d_{1,1}, \ldots, x_n-\frac{1}{2}d_{n,n})
\] 
and 
$H_3\circ H_2\circ H_1\circ G_{(p, A)}\circ h_1=(\varphi_{4,0}, \varphi_{4,1}, \ldots, \varphi_{4,2n})$.    
Then, $\varphi_{4, i}$ $(0\le i\le 2n)$ may be expressed as follows, where 
$d_{i,0}, d_{i, i-n}$ are non-zero 
real numbers obtained in STEP 3 and $\widetilde{d}_i$ are some 
constants.    
\[
\varphi_{4, i}(x_0, x_1, \ldots, x_n)  
= 
\left\{ 
\begin{array}{ll}
x_i^2+\widetilde{d}_i &  (0\le i\le n) \\
 d_{i, 0}x_0+d_{i, i-n}x_{i-n}+\widetilde{d}_i  &  (n+1\le i\le 2n).  
\end{array} 
\right.
\]
%%%%%%%%%%%%%%%%%%%%%%%%%%%%%%%%%%%%%%%%%%%%%%%%%%%%%%%%%%%%%%%%%%% 
\subsection{STEP 5}\label{subsection 2.5} 
%%%%%%%%%%%%%%%%%%%%%%%%%%%%%%%%%%%%%%%%%%%%%%%%%%%%%%%%%%%%%%%%%%% 
The purpose of this step is to delete constant terms by the parallel transformation 
$H_4: \mathbb{R}^{2n+1}\to \mathbb{R}^{2n+1}$ 
defined by 
\[
H_4(X_0, X_1, \ldots, X_{2n})=
(X_0-\widetilde{d}_0, X_1-\widetilde{d}_1, \ldots, X_{2n}-\widetilde{d}_{2n}).
\]    
Set 
$H_4\circ H_3\circ H_2\circ H_1\circ G_{(p, A)}\circ h_1=
(\varphi_{5,0}, \varphi_{5,1}, \ldots, \varphi_{5,2n})$.    
Then, $\varphi_{5, i}$ $(0\le i\le 2n)$ may be expressed as follows. 
%, where $d_{i,j}$ are non-zero 
%constants obtained in STEP 3.    
\[
\varphi_{5, i}(x_0, x_1, \ldots, x_n)  
= 
\left\{ 
\begin{array}{ll}
x_i^2 &  (0\le i\le n) \\
d_{i, 0}x_0+d_{i, i-n}x_{i-n}  &  (n+1\le i\le 2n).  
\end{array} 
\right.
\]
%%%%%%%%%%%%%%%%%%%%%%%%%%%%%%%%%%%%%%%%%%%%%%%%%%%%%%%%%%%%%%%%%%% 
\subsection{STEP 6}\label{subsection 2.6} 
%%%%%%%%%%%%%%%%%%%%%%%%%%%%%%%%%%%%%%%%%%%%%%%%%%%%%%%%%%%%%%%%%%% 
The purpose of this step is to simplify the last $n$ components of 
$H_4\circ H_3\circ H_2\circ H_1\circ G_{(p, A)}\circ h_1$.   
Set 
\[
\widetilde{x}_0=x_0, \widetilde{x}_1=d_{n+1, 0}x_0+d_{n+1, 1}x_{1} , \ldots,   
\widetilde{x}_n=d_{i, 0}x_0+d_{2n, n}x_{n}.
\]       
Then, 
since $d_{i, i-n}\ne 0$ for any $i$ $(n+1\le i\le 2n)$, 
the mapping 
$(x_0, x_1, \ldots, x_n)\mapsto (\widetilde{x}_0, \widetilde{x}_1, \ldots, \widetilde{x}_n)$ 
is a linear transformation of $\mathbb{R}^{n+1}$.   
Thus, setting 
\[
h_2(x_0, x_1, \ldots, x_n)=
\left(x_0, \frac{1}{d_{n+1, 1}}({x}_1-d_{n+1, 0}x_0), \ldots,  \frac{1}{d_{2n, n}}({x}_n-d_{2n, 0}x_0)\right)
\] 
and 
composing $h_2$ to $H_4\circ H_3\circ H_2\circ H_1\circ G_{(p, A)}\circ h_1$, 
the desired form may be obtained as follows:   
\begin{eqnarray*}
{ } & { } & 
H_4\circ H_3\circ H_2\circ H_1\circ G_{(p, A)}\circ h_1\circ h_2(x_0, x_1, \ldots, x_n) \\ 
{ } & = & 
\left(
x_0^2, \frac{1}{d_{n+1,1}^2}(x_1-d_{n+1,0}x_0)^2, \ldots, \frac{1}{d_{2n,n}^2}(x_n-d_{2n,0}x_0)^2, x_1, \ldots, x_n
\right). 
\end{eqnarray*}  
%%%%%%%%%%%%%%%%%%%%%%%%%%%%%%%%%%%%%%%%%%%%%%%%%%%%%%%%%%%%%%%%%%% 
\subsection{STEP 7}\label{subsection 2.7} 
%%%%%%%%%%%%%%%%%%%%%%%%%%%%%%%%%%%%%%%%%%%%%%%%%%%%%%%%%%%%%%%%%%% 
This is the final step.     
Notice that $d_{i,0}\ne 0$ for any $i$ $(n+1\le i\le 2n)$.     
Thus, the following coordinate transformation $H_5$ is well-defined.    
\begin{eqnarray*}
{ } & { } & H_5(X_0, X_1, \ldots, X_{2n}) \\ 
{ } & = & 
\left(
X_0, 
-\frac{d_{n+1, 1}^2}{2d_{n+1,0}}
\left(
X_1 - \frac{d_{n+1,0}^2}{d_{n+1,1}^2}X_0 - \frac{1}{d_{n+1,1}^2}X_{n+1}^2
\right), 
\ldots,  \right. \\ 
{ } & { } & \left. \quad\quad 
-\frac{d_{2n, n}^2}{2d_{2n,0}}
\left(
X_n - \frac{d_{2n,0}^2}{d_{2n,n}^2}X_0 - \frac{1}{d_{2n,n}^2}X_{2n}^2
\right), 
X_{n+1}, \ldots, X_{2n}
\right).   
\end{eqnarray*}
Then, the desired normal form of Whitney umbrella may be obtained as follows:  
\begin{eqnarray*}
{ } & { } & 
H_5\circ H_4\circ H_3\circ H_2\circ H_1\circ G_{(p, A)}\circ h_1\circ h_2(x_0, x_1, \ldots, x_n) \\ 
{ } & = & 
\left(
x_0^2, x_0x_1, \ldots, x_0x_n, x_1, \ldots, x_n
\right). 
\end{eqnarray*}  
\hfill $\Box$
%%%%%%%%%%%%%%%%%%%%%%%%%%%%%%%%%%%%%%%%%%%%%%%%%%
%%%%%%%%%%%%%%%%%%%%%%%%%%%%%%%%%%%%%%%%%%%%%%%%%%    
\section{Proof of the assertion (2) of Theorem \ref{theorem 1}}\label{section 3}
%%%%%%%%%%%%%%%%%%%%%%%%%%%%%%%%%%%%%%%%%%%%%%%%%%
%%%%%%%%%%%%%%%%%%%%%%%%%%%%%%%%%%%%%%%%%%%%%%%%%% 
Set   ${A}_3=(a_{ij})_{0\le i, j\le n-1}$. 
% and ${A}_2=(a_{ij})_{n+1\le i\le 2n, 0\le j\le n}$.     
Taking permutations of coordinates of the target space if necessary,  without loss of generality, 
from the first we may assume that rank$(A_3)=$rank$(A)\le n$. 
%%%%%%%%%%%%%%%%%%%%%%%%%%%%%%%%%%%%%%%%%%%%%%%%%%    
\subsection{STEP 1}\label{subsection 3.1}
%%%%%%%%%%%%%%%%%%%%%%%%%%%%%%%%%%%%%%%%%%%%%%%%%%
The purpose of this step is to delete quadratic terms as many as possible 
by a linear transformation $H_1: \mathbb{R}^{2n+1}\to \mathbb{R}^{2n+1}$ of the following type.    
\begin{eqnarray*}
{ } & { } & H_1(X_0, X_1, \ldots, X_{2n}) \\ 
{ } & = & (X_0, X_1, \ldots, X_{2n}) 
\left(
\begin{array}{cccccccc}
1 & 0 & \cdots & 0 & 0 & \lambda_{0,n} & \cdots & \lambda_{0,2n} \\ 
0 & 1 & \ddots & 0 & 0 & \lambda_{1,n} & \cdots & \lambda_{1,2n} \\
\vdots & \ddots & \ddots & \vdots & \vdots & \vdots & \vdots & \vdots \\  
0 & 0 & \cdots & 1 & 0  & \lambda_{n-2,n} & \cdots & \lambda_{n-2,2n} \\ 
0 & 0 & \cdots & 0 & 1  & \lambda_{n-1,n} & \cdots & \lambda_{n-1,2n} \\ 
%\lambda_{n,0} & \cdots  & \lambda_{n,n} & \lambda_{n,n+1} & \cdots & \lambda_{n,2n} \\ 
0 & 0 & \cdots & 0 & 0 & 1 & 0 & 0 \\ 
0 & 0 & \cdots & 0 & 0 & 0 &\ddots & 0 \\  
%\vdots & \ddots & \vdots & 0 & \ddots & 0 \\ 
0 & 0 & \cdots & 0 & 0 & 0 & 0 & 1  
\end{array}
\right).   
\end{eqnarray*}
For any $i$ $(1\le i\le 2n+1)$, let ${\bf a}_{i-1}$ be the $i$-th row vector of $A$.    
Then, by the assumption of rank$(A_3)=$rank$(A)\le n$, the following claim holds: 
\begin{claim}\label{claim 3.1}
For any $i$ $(n\le i\le 2n+1)$, there exist $\alpha_{0i}, \ldots, \alpha_{(n-1)i}\in \mathbb{R}$ such that 
the following equality holds, 
where ${\bf 0}$ is the $(n+1)$ dimensional zero row vector: 
\[
\sum_{j=0}^{n-1} \alpha_{ji}{\bf a}_j + {\bf a}_i={\bf 0}.   
\] 
\end{claim}
Set $H_1\circ G_{(p, A)}=(\varphi_{1,0}, \varphi_{1,1}, \ldots, \varphi_{1,2n})$.    
Then, $\varphi_{1, i}$ $(0\le i\le 2n)$ may be expressed as follows:    
\[
\varphi_{1, i}(x_0, x_1, \ldots, x_n)  
= 
\left\{ 
\begin{array}{ll}
\sum_{j=0}^n a_{ij}(x_j-p_{ij})^2 &  (0\le i\le n-1) \\
\sum_{j=0}^n b_{ij}x_j + c_i&  (n\le i\le 2n), 
\end{array} 
\right.
\]
where $b_{ij}, c_i$ stands for some constants.   
%%%%%%%%%%%%%%%%%%%%%%%%%%%%%%%%%%%%%%%%%%%%%%%%%%%%%%%%%%%%%%%%%%% 
\subsection{STEP 2}\label{subsection 3.2} 
%%%%%%%%%%%%%%%%%%%%%%%%%%%%%%%%%%%%%%%%%%%%%%%%%%%%%%%%%%%%%%%%%%% 
The purpose of this step is to delete constant terms $c_i$ 
by the parallel transformation $H_2: \mathbb{R}^{2n+1}\to \mathbb{R}^{2n+1}$ 
defined by 
\[
H_2(X_0, X_1, \ldots, X_{2n})=
(X_0, X_1, \ldots, X_{n-1}, X_{n}-c_n, \ldots, X_{2n}-c_{2n}). 
\]    
Set $H_2\circ H_1\circ G_{(p, A)}=(\varphi_{2,0}, \varphi_{2,1}, \ldots, \varphi_{2,2n})$.    
Then, $\varphi_{2, i}$ $(0\le i\le 2n)$ may be expressed as follows.    
\[
\varphi_{2, i}(x_0, x_1, \ldots, x_n)  
= 
\left\{ 
\begin{array}{ll}
\sum_{j=0}^n a_{ij}(x_j-p_{ij})^2 &  (0\le i\le n-1) \\
\sum_{j=0}^n b_{ij}x_j &  (n\le i\le 2n), 
\end{array} 
\right.
\]
%%%%%%%%%%%%%%%%%%%%%%%%%%%%%%%%%%%%%%%%%%%%%%%%%%    
\subsection{STEP 3}\label{subsection 3.3}
%%%%%%%%%%%%%%%%%%%%%%%%%%%%%%%%%%%%%%%%%%%%%%%%%%
The purpose of this step is to construct the bad set $\Sigma_A$. 
Set ${B}=(b_{ij})_{n\le i\le 2n, 0\le j\le n }$.      
Then, $B$ is an $(n+1)\times (n+1)$ matrix.   
\begin{definition}
{\rm 
$\Sigma_{A}$ is the set consisting of points $p\in (\mathbb{R}^{n+1})^{2n+1}$ 
such that $\det B=0$.     
}
\end{definition}
The set $\Sigma_A$ is closed and 
of Lebesgue measure zero since it is an algebraic set.    
Take a point $p\in \left(\mathbb{R}^{n+1}\right)^{2n+1}-\Sigma_A $.    
Then, by the construction of $\Sigma_A$, 
the matrix $B$ has its inverse matrix $B^{-1}$.    
Let $H_3: \mathbb{R}^{2n+1}\to \mathbb{R}^{2n+1}$ be the linear transformation 
defined as follows:       
\begin{eqnarray*}
{ } & { } & H_3(X_0, X_1, \ldots, X_{2n}) \\ 
{ } & = & (X_0, X_1, \ldots, X_{2n}) 
\left(
\begin{array}{c|c}
E_n &  \bigzerou \\ 
\hline 
\bigzerol &  \left(B^T\right)^{-1} 
\end{array}
%\begin{array}{ccccc}
%1 & 0 & \cdots & 0 & {\bf 0} \\ 
%0 & 1 & \ddots & 0 & \vdots \\ 
%0 & 0  & \ddots & 0 & \vdots \\ 
%0 & 0  & \cdots & 1 & {\bf 0} \\
%{\bf 0}^T & {\bf 0}^T & \cdots & {\bf 0}^T & B^{-1} 
%\end{array}
\right), 
\end{eqnarray*}
%where ${\bf 0}^T$  is the $(n+1)$ dimensional zero column vector.     
where $E_n$ is the $n\times n$ unit matrix.   
\par 
Set $H_3\circ H_2\circ H_1\circ G_{(p, A)}=(\varphi_{3,0}, \varphi_{3,1}, \ldots, \varphi_{3,2n})$.
Then, $\varphi_{3, i}$ $(0\le i\le 2n)$ may be expressed as follows.    
\[
\varphi_{3, i}(x_0, x_1, \ldots, x_n)  
= 
\left\{ 
\begin{array}{ll}
\sum_{j=0}^n a_{ij}(x_j-p_{ij})^2 &  (0\le i\le n-1) \\
x_{i-n} &  (n\le i\le 2n).    
\end{array} 
\right.
\]
%%%%%%%%%%%%%%%%%%%%%%%%%%%%%%%%%%%%%%%%%%%%%%%%%%%%%%%%%%%%%%%%%%% 
\subsection{STEP 4}\label{subsection 3.4} 
%%%%%%%%%%%%%%%%%%%%%%%%%%%%%%%%%%%%%%%%%%%%%%%%%%%%%%%%%%%%%%%%%%% 
The purpose of this step is to delete remaining constant terms by the parallel transformation $H_4: \mathbb{R}^{2n+1}\to \mathbb{R}^{2n+1}$ 
defined by 
\begin{eqnarray*}
{ } & { } & H_4(X_0, X_1, \ldots, X_{2n}) \\ 
{ } & = & (X_0-\sum_{j=0}^n a_{0j}p_{0j}^2, X_1-\sum_{j=0}^n a_{1j}p_{1j}^2, \ldots, 
X_{n-1}-\sum_{j=0}^n a_{(n-1)j}p_{(n-1)j}^2,  \\ 
{ } & { } & \qquad \qquad X_n, X_{n+1}\ldots, X_{2n}). 
\end{eqnarray*}
Set 
$H_4\circ H_3\circ H_2\circ H_1\circ G_{(p, A)}=
(\varphi_{4,0}, \varphi_{4.1}, \ldots, \varphi_{4,2n})$.    
Then, $\varphi_{4, i}$ $(0\le i\le 2n)$ may be expressed as follows.  
% where, $d_{i,j}$ are non-zero 
%constants obtained in STEP 3.    
\[
\varphi_{4, i}(x_0, x_1, \ldots, x_n)  
= 
\left\{ 
\begin{array}{ll}
\sum_{j=0}^n a_{ij}(x_j^2-2p_{ij}x_j) &  (0\le i\le n-1) \\
x_{i-n} &  (n\le i\le 2n).    
\end{array} 
\right.
\]
%%%%%%%%%%%%%%%%%%%%%%%%%%%%%%%%%%%%%%%%%%%%%%%%%%%%%%%%%%%%%%%%%%% 
\subsection{STEP 5}\label{subsection 3.5} 
%%%%%%%%%%%%%%%%%%%%%%%%%%%%%%%%%%%%%%%%%%%%%%%%%%%%%%%%%%%%%%%%%%% 
This is the final step.     
Let $H_5: \mathbb{R}^{2n+1}\to \mathbb{R}^{2n+1}$ be the coordinate transformation defined by 
\begin{eqnarray*}
{ } & { } & H_5(X_0, X_1, \ldots, X_{2n}) \\ 
{ } & = & 
\left(
X_n, X_{n+1}, \ldots, X_{2n}, X_0-\sum_{j=0}^n a_{0j}\left(X_{n+j}^2-2p_{0j}X_{n+j}\right), 
 \right.
\\ 
{ } & { } &  
\qquad 
X_1-\sum_{j=0}^n a_{1j}\left(X_{n+j}^2-2p_{1j}X_{n+j}\right),  
\ldots,  \\ 
{ } & { } & \left.
\qquad\qquad 
X_{n-1}-\sum_{j=0}^n a_{(n-1)j}\left(X_{n+j}^2-2p_{(n-1)j}X_{n+j}\right) 
\right).  
\end{eqnarray*}
Then, we have the following:     
\[
H_5\circ H_4\circ H_3\circ H_2\circ H_1\circ G_{(p, A)}(x_0, x_1, \ldots, x_n) = 
\left(
x_0, \ldots, x_n, 0, \ldots, 0
\right).   
\]
%Finally, composing the permutation of the target space with 
%$H_3\circ H_2\circ H_1\circ G_{(p, A)}$ yields the desired inclusion 
%$(x_0, \ldots, x_n)\mapsto (x_0, \ldots, x_n, 0, \ldots, 0)$.    
\hfill $\Box$
%%%%%%%%%%%%%%%%%%%%%%%%%%%%%%%%%%%%%%%%%%%%%%%%%%%%%%%%%%%%%%%%%%% 
\subsection{REMARK}\label{subsection 3.6} 
%%%%%%%%%%%%%%%%%%%%%%%%%%%%%%%%%%%%%%%%%%%%%%%%%%%%%%%%%%%%%%%%%%% 
As a by-product of the proof of the assertion (2) of 
Theorem \ref{theorem 1} given in this section,  
%Section \ref{section 2} and Section \ref{section 3}, 
we have the following:   
\begin{theorem}
\label{theorem 3}
Let $n$ be a positive integer and $k$ be an integer such that $k>2n$.    
Let $A=(a_{ij})_{0\le i\le k, 0\le j\le n}$ be a $(k+1)\times (n+1)$ matrix with non-zero entries.   
Then, there exists a closed subset $\Sigma_A\subset (\mathbb{R}^{n+1})^{2n+1}$ 
with Lebesgue measure zero 
such that for any $p=(p_0, \ldots, p_{2n})\in  (\mathbb{R}^{n+1})^{2n+1}-\Sigma_A$, 
$G_{(p, A)}$ is $\mathcal{A}$-equivalent to the inclusion  
$(x_0, \ldots, x_n)\mapsto (x_0, \ldots, x_n, 0, \ldots, 0).$
\end{theorem}   
%%%%%%%%%%%%%%%%%%%%%%%%%%%%%%%%%%%%%%%%%%%%%%%%%% 
%%%%%%%%%%%%%%%%%%%%%%%%%%%%%%%%%%%%%%%%%%%%%%%%%% 
%%%%%%%%%%%%%%%%%%%%%%%%%%%%%%%%%%%%%%%%%%%%%%%%%%    
\section{Proof of Theorem \ref{theorem 2}}\label{section 4}
%%%%%%%%%%%%%%%%%%%%%%%%%%%%%%%%%%%%%%%%%%%%%%%%%%
%%%%%%%%%%%%%%%%%%%%%%%%%%%%%%%%%%%%%%%%%%%%%%%%%%  
It is sufficient to show the following proposition.   
\begin{proposition}\label{proposition 3}
Let $X$ be the set consisting of $(2n+1)\times (n+1)$ matrices $A$ with non-zero entries such that 
the rank of $A$ is $n+1$.      
Then, for any $A\in X$ there exist an open neighborhood $U_A$ of $A$ in $X$ 
and an open set 
$V\subset (\mathbb{R}^{n+1})^{2n+1}$ such that 
for any $p\in V$ there exists $B\in U_A$ satisfying that $G_{(p, B)}$ is an unstable mapping.    
\end{proposition}
\underline{Proof of Proposition \ref{proposition 3}.}\quad 
Let $A=(a_{ij})_{0\le i\le 2n, 0\le j\le n}$ be a given element of $X$.     
Taking a permutation of coordinates of the target space if necessary, without loss of generality, 
from the first we may assume that the rank of $A_1$ is $n+1$, where $A_1=(a_{ij})_{0\le i, j\le n}$ as defined in 
the beginning of Section \ref{section 2}.       
\par 
Let $\widetilde{A}=(\widetilde{a}_{ij})_{0\le i\le 2n, 0\le j\le n}\in X$ be a sufficiently near 
matrix to $A$ such that $\widetilde{a}_{ij}=a_{ij}$ if $0\le i, j \le n$.    Set 
$\widetilde{A}_2=(\widetilde{a}_{ij})_{n+1\le i\le 2n, 0\le j\le n}$.   
For the matrix $\widetilde{A}$, consider the 
matrix $\Lambda_2=(\lambda_{i,j})_{0\le i\le n, n+1\le j\le 2n}$ defined in STEP 1 of Section \ref{section 2}:   
\[
\Lambda_2=-(A_1^T)^{-1}\widetilde{A}_2^T.   
\]
%where $A_2=(a_{ij})_{n+1\le i\le 2n, 0\le j\le n}$ is as in Section \ref{section 2}.   
Moreover, as in Section \ref{section 2}, for any $i,j$ $(n+1\le i\le 2n, 0\le j\le n)$, 
consider $b_{ij}$ for the matrix $\widetilde{A}$.    
\[
b_{ij}=-2\left(\sum_{k=0}^n \lambda_{k,i}a_{kj}p_{kj}+\widetilde{a}_{ij}p_{ij}\right),  
\]   
where the real number $\lambda_{k, i}$ is the $(k, i-n)$ component of $\Lambda_2$.    
Notice that $\lambda_{k,i}$ in $b_{ij}$ is a linear function with 
respect to $\widetilde{a}_{ij}$ $(n+1\le i\le 2n, 0\le j\le n)$.   
% 
%$\widetilde{a}_{ij}$ $(n+1\le i\le 2n, 0\le j\le n)$.      
Thus, for any $i_0, j_0$ $(n+1\le i_0\le 2n, 0\le j_0\le n)$, the following function $\widetilde{\psi}_{i_0,j_0}$ is a rational function 
with variables 
%$\widetilde{a}_{ij}$ $(n+1\le i\le 2n, 0\le j\le n)$ and 
$p_{00}, \ldots, p_{nn}, \widetilde{a}_{(n+1)0}, \ldots, \widetilde{a}_{(2n)n}$.  
%  $(0\le k, j \le n)$.    
\[
\widetilde{\psi}_{i_0,j_0}(p_{00}, \ldots, p_{nn}, \widetilde{a}_{(n+1)0}, \ldots, \widetilde{a}_{(2n)n})  
 =  \frac{-\sum_{k=0}^n \lambda_{k,i_0}a_{kj_0}p_{kj_0}}{\widetilde{a}_{i_0j_0}}.   
\]
\par 
We would like to show that there exist a matrix $\widetilde{A}$ of the above type 
which is sufficiently near $A$ and an open set $V$ 
in $(\mathbb{R}^{n+1})^{2n+1}$ such that $b_{ij}=0$ for any point $p\in V$.       
In order to do so, we consider the mapping 
\[\Psi=(\psi_{i_0,j_0})_{0\le i_0\le 2n, 0\le j_0\le n}: 
\mathbb{R}^{(n+1)^2}\times (\mathbb{R}-\{{\bf 0}\})^{n(n+1)}\to \mathbb{R}^{(2n+1)(n+1)}
\] 
defined as follows:   
\begin{eqnarray*}
{ } & { } & \psi_{i_0,j_0}(q_{00}, \ldots, q_{nn}, c_{(n+1)0}, \ldots, c_{(2n)n} ) \\   
{ } & = &  
\left\{ 
\begin{array}{ll}
q_{i_0j_0}   & (0\le i_0\le n) \\
\widetilde{\psi}_{i_0,j_0}(q, c)  &  (n+1\le i_0\le 2n),    
\end{array} 
\right.
\end{eqnarray*}
where $q=(q_{00}, \ldots, q_{nn})$ and $c=(c_{(n+1)0}, \ldots, c_{(2n)n})$.     
\begin{definition}
{\rm 
\[
\widetilde{\Sigma}=\left\{\left.(q,c)\in \mathbb{R}^{(n+1)^2}\times (\mathbb{R}-\{{\bf 0}\})^{n(n+1)}\; \right|\; \det J\Psi(q, c)=0\right\}, 
\]
where $J\Psi(q, c)$ is the Jacobian matrix of $\Psi$ at $(q,c)$.   
}
\end{definition}
Since $\psi_{i_0,j_0}$ is a rational function, $\widetilde{\Sigma}$ is 
a semi-algebraic subset of Lebesgue measure zero.    
Let $(q_0, c_0)$ be a point outside $\widetilde{\Sigma}$.     
Since $\widetilde{\Sigma}$ is closed and of Lebesgue measure zero, 
we may assume that $c_0$ is sufficiently near $A_2$.    
Since $\det J\Psi(q_0, c_0)\ne 0$, by the inverse function theorem, 
there exist open neighborhoods $U_1\subset\mathbb{R}^{(n+1)^2} $ of $q_0$, 
$U_2\subset (\mathbb{R}-\{{\bf 0}\})^{n(n+1)} $ of $c_0$ and 
$V\subset \mathbb{R}^{(2n+1)(n+1)}$ of $\Psi(q_0, c_0)$ such that the restriction 
$\Psi|_{U_1\times U_2}: U_1\times U_2\to V$ is a $C^\infty$ diffeomorphism.     
In particular, we have the following:   
\begin{claim}\label{claim 4.1}
For any $p=(p_{00}, \ldots, p_{(2n)n})\in V$ there exists 
a matrix $\widetilde{A}_2\in U_2$ such that 
\[
p_{i_0j_0}=\widetilde{\psi}_{i_0,j_0}(p_{00}, \ldots, p_{nn}, \widetilde{A}_2  %\widetilde{a}_{(n+1)0}, \ldots, a_{(2n)n}
)
\]
for any $i_0, j_0$ $(n+1\le i_0\le 2n, 0\le j_0\le n)$.   
\end{claim}
Claim \ref{claim 4.1} implies the following:   
\begin{claim}\label{claim 4.2}
For any $p=(p_{00}, \ldots, p_{(2n)n})\in V$ there exists 
a matrix $\widetilde{A}_2\in U_2$ such that $b_{i_0j_0}=0$ 
for any $i_0, j_0$ $(n+1\le i_0\le 2n, 0\le j_0\le n)$.   
\end{claim}
Claim 4.2 shows that the image of $H_1\circ G_{(p, \widetilde{A})}$ must be inside $\mathbb{R}^{n+1}\times \{{\bf 0}\}$ 
where $H_1$ is the linear transformation given in STEP 1 of Section \ref{section 2}.    
Thus, by the classification of stable singularity for map-germ $\mathbb{R}^{n+1}\to \mathbb{R}^{2n+1}$ due to 
Whitney (\cite{whitney}), if $G_{(p, \widetilde{A})}$ is stable, then it must be $\mathcal{A}$-equivalent to the inclusion 
$(x_0, \ldots, x_n)\mapsto (x_0, \ldots, x_n, 0, \ldots, 0)$.      In particular, $G_{(p, \widetilde{A})}$ must be non-singular.   
However, it is easily seen that $G_{(p, \widetilde{A})}$ is singular.  Hence, $G_{(p, \widetilde{A})}$ is unstable.   
\hfill $\Box$
%%%%%%%%%%%%%%%%%%%%%%%%%%%%%%%%%%%%%%%%%%%%%%%%%%%%%%%%%%%%%%%%%%%
%%%%%%%%%%%%%%%%%%%%%%%%%%%%%%%%%%%%%%%%%%%%%%%%%%%%%%%%%%%%%%%%%%% 
\subsection{REMARK}\label{subsection 4.1} 
%%%%%%%%%%%%%%%%%%%%%%%%%%%%%%%%%%%%%%%%%%%%%%%%%%%%%%%%%%%%%%%%%%%
%%%%%%%%%%%%%%%%%%%%%%%%%%%%%%%%%%%%%%%%%%%%%%%%%%%%%%%%%%%%%%%%%%% 
\begin{enumerate}
\item 
As a by-product of the proof of Theorem \ref{theorem 2} given in this section,  
%Section \ref{section 2} and Section \ref{section 3}, 
we have the following:   
\begin{theorem}
\label{theorem 4}
Let $n$ be a positive integer and $k$ be an integer such that $k>2n$.   
Then, there does not  exist a semi-algebraic subset $\Sigma\subset (\mathbb{R}^{n+1})^{k+1}$  with Lebesgue measure zero 
such that for any point $p=(p_0, p_1, \ldots, p_{k})\in  (\mathbb{R}^{n+1})^{k+1}-\Sigma$ and any $(k+1)\times (n+1)$ matrix $A$ with non-zero entries, 
$G_{(p, A)}$ is $\mathcal{A}$-equivalent to 
%the following mapping:   
the inclusion 
$(x_0, \ldots, x_n)\mapsto (x_0, \ldots, x_n, 0, \ldots, 0).$
\end{theorem} 
\item The proof of Theorem \ref{theorem 2} has one more advantage.   
It makes clear the reason why 
we can expect the existence of a universal bad set $\Sigma$ in the case $n=k=1$.  
\end{enumerate}
%%%%%%%%%%%%%%%%%%%%%%%%%%%%%%%%%%%%%%%%%%%%%%%%%%%%%%%    
\section*{Acknowledgements}
The authors are grateful to M.~A.~S.~Ruas and R.~Oset Sinha for their kind suggestions and advices. 
T.~Nishimura is partially supported by JSPS-CAPES under the JAPAN-BRAZIL research cooperative program.   
%R.~Oset Sinha is partially supported by FAPESP
%grant no. 2013/02381-1 and DGCYT and FEDER grant no. MTM2012-33073.
%M. A. S.~Ruas is partially supported by CNPq grant no.
%305651/2011-0.

%%%%%%%%%%%%%%%%%%%%%%%%%%%%%%%%%%%%%%%%%%%%%%%%%%   
%%%%%%%%%%%%%%%%%%%%%%%%%%%%%%%%%%%%%%%%%%%%%%%%%%  

\end{document}